\documentclass[12pt]{amsart}
\usepackage{amssymb, amsthm, enumerate}
\usepackage{graphicx}
\usepackage{indentfirst}
\newtheorem{Def}{Definition}
\newtheorem{Lem}{Lemma}

\newtheorem{Thm}{Theorem}
\newtheorem{Cor}{Corollary}
\newtheorem{Rem}{Remark}
\newenvironment{Pf}{ Proof.}{\(\square\)}
\newtheorem{Exm}{Example}
\newtheorem{Exc}{Exercise}

\begin{document}
\author{Csaba Vincze}
\footnotetext[2]{{\bf MR subject classification:} primary 51B20,
secondary 52A40}
\footnotetext[3]{Cs. Vincze is supported by EFOP 3.6.2-16-2017-00015.}
\address{Institute of Mathematics, University of Debrecen,\\ P.O.Box 12, H-4010 Debrecen, Hungary}
\email{csvincze@science.unideb.hu}
\title[Lazy orbits: an optimization problem on the sphere...]{Lazy orbits: an optimization problem on the sphere}
\begin{abstract}
Non-transitive subgroups of the orthogonal group play an important role in the non-Euclidean geometry. If $G$ is a closed subgroup in the orthogonal group such that the orbit of a single Euclidean unit vector does not cover the (Euclidean) unit sphere centered at the origin then there always exists a non-Euclidean Minkowski functional such that the elements of $G$ preserve the Minkowskian length of vectors. In other words the Minkowski geometry is an alternative of the Euclidean geometry for the subgroup $G$. It is rich of isometries if $G$ is "close enough" to the orthogonal group or at least to one of its transitive subgroups. The measure of non-transitivity is related to the Hausdorff distances of the orbits under the elements of $G$ to the Euclidean sphere. Its maximum/minimum belongs to the so-called lazy/busy orbits, i.e. they are the solutions of an optimization problem on the Euclidean sphere. The extremal distances allow us to characterize the reducible/irreducible subgroups. We also formulate an upper and a lower bound for the ratio of the extremal distances. 

As another application of the analytic tools we introduce the rank of a closed non-transitive group $G$. We shall see that if $G$ is of maximal rank then it is finite or reducible. Since the reducible and the finite subgroups form two natural prototypes of non-transitive subgroups, the rank seems to be a fundamental notion in their characterization. Closed, non-transitive groups of rank $n-1$ will be also characterized. Using the general results we classify all their possible types in lower dimensional cases $n=2, 3$ and $4$. 

Finally we present some applications of the results to the holonomy group of a metric linear connection on a connected Riemannian manifold.
\end{abstract}
\maketitle

\section{Introduction}

In the preamble to his fourth problem presented at the International Mathematical Congress in Paris (1900) Hilbert suggested the examination of geometries standing next to Euclidean one in the sense that they satisfy much of Euclidean's axioms except some (tipically one) of them. In the classical non-Euclidean geometry the axiom taking to fail is the fameous parallel postulate. Another type of geometry standing next to Euclidean one is the geometry of normed spaces or, in a more general context, the geometry of Minkowski spaces \cite{PaivaThompson}. The crucial test is not the parallelism but the congruence via the group of linear isometries. Consider the Euclidean coordinate space $\mathbb{R}^n$ and let $G\subset O(n)$ be a closed non-transitive subgroup in the orthogonal group. We are going to construct a compact convex body $K$ containing the origin in its interior such that

\begin{itemize}
\item[(K1)] $K$ is not a unit ball with respect to any inner product (ellipsoid problem),

\item[(K2)] $K$ is invariant under the elements of the subgroup $G$,

\item[(K3)] its boundary $\partial K$ is a smooth hypersurface (regularity condition).
\end{itemize}

The Minkowski functional induced by $K$ as a unit ball makes the vector space a non-Euclidean Minkowski space in the sense of condition (K1). The second condition (K2) says that $G$ is a subgroup of the linear isometry group with respect to the Minkowski functional. In other words the Minkowski geometry is the alternative of the Euclidean geometry for the subgroup $G$. The regularity condition (K3) allows us to apply the standard differential geometric tools in this situation: the central role of the theory is played by the Hessian of the Minkowskian norm square \cite{BCS}. In case of a differentiable manifold equipped with a Riemannian metric, the subgroup $G$ can be be interpreted as the holonomy group of a metric linear connection and the alternative geometry is the Finsler geometry: \emph{instead of the Euclidean spheres in the tangent spaces, the unit vectors form the boundary of general convex sets containing the origin in their interiors} (M. Berger). 

One of the main applications of generalized conics' theory is to provide a method for the construction of convex bodies satisfying conditions (K1) - (K3); \cite{VA_1}, see also \cite{VA_2}. In a more general sense a conic is a subset of the space all of whose points have the same average Euclidean distance from the elements of the focal set $D$, i.e. the function
\begin{equation}
v \mapsto f_D(v):=\int_{D} u\circ d_2(v,w)\, d_{\mu} w
\end{equation}
is constant, where $d_2$ means the Euclidean distance function on $\mathbb{R}^n$, $D\subset \mathbb{R}^n$ is a compact subset with a finite positive measure with respect to $\mu$ and $u\colon \mathbb{R}\to \mathbb{R}$ is a strictly monote increasing convex function satisfying the initial condition $u(0)=0$. Depending on the context, a generalized conic  means the set all of whose points satisfy
$$f_D(v)=c\ \ \textrm{\ or\ }\ \ f_D(v)\leq c.$$
It can be proved that $f_D$ is a convex and, consequently, continuous function on the coordinate space satisfying the growth condition
$$\liminf_{\|v\|_2\to \infty}\ \frac{f_D(v)}{\|v\|_2}>0.$$
Therefore any sublevel subset $f_D(v)\leq c$ is a convex compact subset in $\mathbb{R}^n$. 

\begin{Rem}{\emph{Using the distance coming from the taxicab norm instead of the Euclidean one, we have another application in geometric tomography, see \cite{VA_3}. In terms of the probability theory, the value of the function measuring the average distance at $v$ is the expectable value of the random variable $d(v, \xi)$ for any uniformly distributed random (vector) variable $\xi$ on the focal set $D$; note that the divison by the measure of the integration domain has no any influence on the level sets. For different distributions see \cite{VNBN}.}}
\end{Rem}

To satisfy (K2) it is natural to choose a $G$-invariant integration domain $D$ equipped with a $G$-invariant measure $\mu$. The idea is illustrated by the following theorem due to L. Bieberbach.

\begin{Thm} \emph{\cite{C}} The holonomy group of any flat compact Riemannian manifold is finite.
\end{Thm}

Let $M$ be a flat compact Riemannian manifold and choose a point $p\in M$. In the sense of Bieberbach's theorem we can find a finite system $w_1, \ldots, w_m$ of elements in the tangent space $T_pM$ which is invariant under the holonomy group $G:=H_p$ of the L\'{e}vi-Civita connection. Choosing the counting measure on the set of the elements $w_i$'s ($i=1, \ldots, m$) and $u(t):=t$ it can be easily seen that 
$$K_p(w_1, \ldots, w_m):=\{v\in T_pM\ | \ \sum_{i=1}^m d_2(v,w_i)\leq c \}\subset T_pM$$
satisfies conditions (K1)-(K3) provided that the constant $c$ is large enough for $K_p(w_1, \ldots, w_m)$ to contain $w_1, \ldots, w_m$ in its interior (see the regularity condition). The conic $K_p(w_1, \ldots, w_m)$ is called a \emph{polyellipsoid} with focuses  $w_1, \ldots, w_m$. In case of $m=2$ we get back the classical notion of ellipses in the Euclidean plane. Since $K_p(w_1, \ldots, w_m)$ is invariant under the elements of the holonomy group $H_p$ we have an induced Minkowski functional in the tangent space $T_pM$ such that the elements of $H_p$ are Minkowskian isometries. Especially $H_p$ is the holonomy group of the L\'{e}vi-Civita connection and the $H_p$-invariance allows us to extend $K_p(w_1, \ldots, w_m)$ by parallel transports to any point of the manifold; see Z. I. Szab\'{o}'s idea \cite{Szab1}. The smoothly varying family of compact convex bodies provides a non-Riemannian metric environment for the L\'{e}vi-Civita connection: the Minkowski functionals induced by the polyellipsoids in the tangent spaces constitute a so-called Finslerian fundamental function such that the parallel transports with respect to the L\'{e}vi-Civita connection preserve the Finslerian length of tangent vectors. Especially the space is a so-called \emph{locally Minkowski manifold} because of the flatness of the connection. 

In general the holonomy group of a metric linear connection is not finite. To adopt the previous method to the general situation we should develop the theory of conics with infinitely many focal points by using integration instead of finite sums; see \cite{VA_2}. The main result says that if we have a closed subgroup $G$ in the orthogonal group such that the orbit of a single Euclidean unit vector does not cover the (Euclidean) unit sphere centered at the origin then there always exists a non-Euclidean Minkowski functional induced by a generalized conic such that the elements of $G$ preserve the Minkowskian length of vectors. In other words the Minkowski geometry is an alternative of the Euclidean geometry for the subgroup $G$. It is rich of isometries if the group $G$ is "close enough" to the orthogonal group $O(n)$ or at least to one of its transitive subgroups. The following list \cite{GrayGreen} shows the compact connected Lie subgroups\footnote{The Closed Subgroup Theorem \cite{Lee} states that $G$ (as a topologically closed sub\-group) is a Lie subgroup. Since the orthogonal group is compact, so is $G$ with at most finitely many components; for the unit components see the list above.} which are transitive on the Euclidean unit sphere $S_{n-1}\subset \mathbb{R}^n$. In case of these groups there are no alternatives of the Euclidean geometry. For the classification see \cite{MS},  \cite{Borel1} and \cite{Borel2}. 

\vspace{0.3cm}
\noindent
\begin{center}
\small{
\begin{tabular}{|l|l|l|l|l|l|}
\hline
&&&&&\\
SO(n)&SO(n)&SO(n)&SO(7)&SO(8)&SO(16)\\
\hline
&&&&&\\
$-$&U(2k+1)&U(2k)&$-$&U(4)&U(8)\\
\hline
&&&&&\\
$-$&SU(2k+1)&SU(2k)&$-$&SU(4)&SU(8)\\
\hline
&&&&&\\
$-$&$-$&Sp(k)&$-$&Sp(2)&Sp(4)\\
\hline
&&&&&\\
$-$&$-$&Sp(k)$\cdot$SO(2)&$-$&Sp(2)$\cdot$SO(2)&Sp(4)$\cdot$SO(2)\\
\hline
&&&&&\\
$-$&$-$&Sp(k)$\cdot$Sp(1)&$-$&Sp(2)$\cdot$Sp(1)&Sp(4)$\cdot$Sp(1)\\
\hline
&&&&&\\
$-$&$-$&$-$&G$_2$&Spin(7)&Spin(9)\\
\hline
&&&&&\\
n=2k+1$\neq$ 7&n=2(2k+1)&n=4k$\neq$ 8, 16&n=7&n=8&n=16\\
\hline
\end{tabular}}
\end{center}

\vspace{0.3cm}
The paper is devoted to the investigation of non-transitive subgroups. The measure of non-transitivity is introduced by the Hausdorff distances of the orbits under the elements of $G$ to the Euclidean sphere. Its maximum/minimum belongs to the so-called lazy/busy orbits, i.e. they are the solutions of an optimization problem on the Euclidean sphere. The extremal distances allow us to characterize the reducible/irreducible subgroups (Theorem 2). We also formulate an upper and a lower bound for the ratio of the extremal distances (Theorem 3). 

As another application of the analytic tools we introduce the rank of a group $G$. The rank provides us a quadratic upper bound for the dimension of $G$ as a Lie subgroup of the orthogonal group (Theorem 7). We shall see that if $G$ is of maximal rank then it is finite or reducible (Theorem 4). Since the reducible and the finite subgroups form two natural prototypes of non-transitive subgroups, the rank seems to be a fundamental notion in their characterization. Closed, non-transitive groups of rank $n-1$ will be also characterized (Theorem 8). Using the general results we classify all their possible types in lower dimensional cases $n=2, 3$ and $4$ (subsections 5.1 - 5.5). 

Finally we present some applications of the results to the holonomy group of a metric linear connection $\nabla$ on a connected Riemannian manifold. The alternative of the Riemannian geometry for $\nabla$ is called Finsler geometry (see the theory of generalized Berwald manifolds \cite{H1}, \cite{Szab2}, \cite{V10}, \cite{Vin11} and \cite{Wag1}).  

\section{An optimization problem: minimax points of orbits on the Euclidean unit sphere} Let $\mathbb{R}^n$ be the standard Euclidean space equipped with the canonical inner product, norm and distance function:
$$\langle v,w \rangle=v^1w^1+\ldots+v^nw^n,\ \ \|v\|:=\sqrt{(v^1)^2+\ldots+(v^n)^2},\ \ d(v,w)=\|v-w\|.$$
Let
$$S_{n-1}:=\{ \| v\|=1\ | \ v\in \mathbb{R}^n\}$$
be the Euclidean unit sphere centered at the origin and consider a subgroup $G\subset O(n)$. For any Euclidean unit vector $v\in S_{n-1}$ let us define the function
$$f_v \colon \mathbb{R}^n\to \mathbb{R},\ \ f_v(w):=\sup_{g\in G} \| w-g(v)\|.$$

\begin{Lem} For any $v\in S_{n-1}$, the mapping $f_v$ is convex and the Lipschitz property 
\begin{equation}
\label{lip}
\|f_v(w)-f_v(z)\|\leq \|w-z\|
\end{equation}
is satisfied.
\end{Lem}

\begin{Pf}
The convexity follows from the convexity of the norm and the sup - functions:
$$\|\lambda w+(1-\lambda)z-g(v)\|\leq \lambda \|w-g(v)\|+(1-\lambda)\|z-g(v)\|,$$
where $0\leq \lambda \leq 1$. On the other hand
$$\|w-g(v)\|\leq \|w-z\|+\|z-g(v)\|\ \ \Rightarrow\ \ f_v(w)\leq \|w-z\|+f_v(z).$$
Changing the role of $w$ and $z$ the proof is finished. 
\end{Pf}

\begin{Cor}
For any $v\in S_{n-1}$, the mapping $f_v$ takes both its minimum and maximum on the unit sphere.
\end{Cor}

\begin{Pf}
It is well-known that the convexity implies the continuity at any interior point of the domain.
\end{Pf}

\begin{Def}
\label{minimaxpoint}
The minimax point of the orbit $P_G(v):=\{ g(v)\ | \ g\in G\}$
of the unit vector $v$ under the subgroup $G$ is the solution of the optimization problem
\begin{equation}
\label{op1}
\inf f_v\ \ \textrm{subject to}\ \ \|w\|=1;
\end{equation}
if $w$ is a minimax point then  $m_v:=f_v(w)$ is called the minimax value of $P_G(v)$. 
\end{Def}

Some monotonicity and stability properties can be summarized as follows:
\begin{itemize}
\item[(m1)] if $G \subset \tilde{G}$ then $f_v\leq \tilde{f}_v$ and $f_v=\bar{f}_v$, where  $\bar{f}_v$ belongs to the topological closure $\bar{G}$ of the subgroup $G$,
\item[(m2)] for any unit element $v\in S_{n-1}$ the set of the minimax points is a closed, and consequently, compact subset of $S_{n-1}$ as the intersection of $f_v^{-1}(m_v)$ and $S_{n-1}$,
\item[(m3)] the set of the minimax points is invariant under the action of the topological closure $\bar{G}$ of the subgroup $G$,
\item[(m4)] $m_v=m_{g(v)}$ for any $g\in \bar{G}$.
\item[(m5)] if $w$ is the minimax point of the orbit $P_G(v)$ then $m_w\leq m_v$ because
$$m_w=\inf_{z\in S_{n-1}}f_w(z)\leq f_w(v)=\sup_{g\in G} \| v-g(w)\|= \sup_{g\in G} \| g^{-1}(v)-w\|=f_v(w)=m_v,$$
where $g$ runs through the elements of $G$. So does $g^{-1}$.
\item[(m6)] $m_v=m_w$ if and only if $w$ is the minimax point of $P_G(v)$ and vice versa. 
\end{itemize}

It can be easily seen that (\ref{op1}) is equivalent to the optimization problem
\begin{equation}
\label{op2}
\inf f_v^2=\inf \left(\sup_{g\in G} \| w-g(v)\| \right)^2=\inf \sup_{g\in G} \| w-g(v)\|^2\ \ \textrm{subject to}\ \ \|w\|=1,
\end{equation}
where
$$f_v^2(w)=\sup_{g\in G} \| w-g(v)\|^2$$
\begin{Lem}
Let $c\colon I\to \mathbb{R}^n$ be a continuously differentiable space curve defined on the interval $I$ containing the origin in its interior. If $c(0)=w$ and $c'(0)=z$ then
\begin{equation}
\label{derminmax}
\lim_{t\to 0^+}\frac{f_v^2(c(t))-f_v^2(c(0))}{t}=\lim_{t\to 0^+}\frac{f_v^2(w+tz)-f_v^2(w)}{t}=:(f_v^2)'(w;z),
\end{equation}
where $(f_v^2)'(w,z)$ is the one-sided directional derivative of the function $f_v^2$ at $w$ in direction $z$.
\end{Lem}

\begin{Pf} First of all observe that for any $t>0$
$$\left |\frac{f_v^2(c(t))-f_v^2(w+tz)}{t}\right|=\left(f_v(c(t))+f_v(w+tz)\right)\left | \frac{f_v(c(t))-f_v(w+tz)}{t}\right |\leq$$
$$ \left(f_v(c(t))+f_v(w+tz)\right) \left |z-\frac{c(t)-w}{t} \right |=\left(f_v(c(t))+f_v(w+tz)\right) \left |z-\frac{c(t)-c(0)}{t} \right |$$
because of the Lipschitz property (\ref{lip}). Therefore 
$$\left |\frac{f_v^2(c(t))-f_v^2(w+tz)}{t}\right|\to 0$$
as $t\to 0^+$ and, consequently,
$$\lim_{t\to 0^+}\frac{f_v^2(c(t))-f_v^2(c(0))}{t}=\lim_{t\to 0^+}\frac{f_v^2(c(t))-f_v^2(w+tz)}{t}+\lim_{t\to 0^+}\frac{f_v^2(w+tz)-f_v^2(w)}{t}=$$
$$\lim_{t\to 0^+}\frac{f_v^2(w+tz)-f_v^2(w)}{t}=(f_v^2)'(w,z)$$
as was to be proved; note that the one-sided directional derivative $(f_v^2)'(w,z)$ exists because of the convexity of the function. 
\end{Pf}

\begin{Cor}
\label{minimax2}
If $w$ is the minimax point of the orbit $P_G(v)$ then $(f_v^2)'(w; z)\geq 0$ for any tangent vector $z$ to $S_{n-1}$ at $w$.
\end{Cor}

Let $w$ be the minimax point of the orbit $P_G(v)$; the set 
\begin{equation}
\label{basicset}
T_f(v,w):=\{g\in \bar{G}\ | \ f_v(w)=\| w-g(v)\| \}
\end{equation}
plays the central role in the derivative of the sup-function \cite{Kawasaki}. $T_f(v,w)$ contains the group elements in $\bar{G}$ sending $v$ to the elements of the orbit where the minimax value is attained at. We have 
\begin{equation}
\label{supderivative}
(f_v^2)'(w,z)=\max_{g\in T_f(v,w)} {\| {\bf x}-g(v)\|^2}'(w; z)=2\max_{g\in T_f(v,w)} \langle w-g(v),z \rangle=
\end{equation}
$$2\max_{g\in T_f(v,w)} -\langle g(v)-w,z \rangle.$$

\begin{Cor}
\label{minimax3}
Let $w$ be the minimax point of the orbit $P_G(v)$; for any tangent vector $z$ to $S_{n-1}$ at $w$
\begin{equation}
\label{minimax4}
\max_{g\in T_f(v,w)} - \langle  g(v)-w, z\rangle \geq 0\ \ \Rightarrow\ \ \min_{g\in T_f(v,w)} \langle g(v)-w, z\rangle \leq 0.
\end{equation}
\end{Cor}

Note that $\langle w, z \rangle=0$ because $z$ is tangential to $S_{n-1}$ at $w$. Inequality (\ref{minimax4}) implies that for any tangent vector $z\in T_w S_{n-1}$ there is an element $g\in T_f(v,w)$ such that $z$ and $g(v)-w$ encloses an angle greater or equal than $90^{\circ}$. 

\begin{figure}
\centering
\includegraphics[scale=0.45]{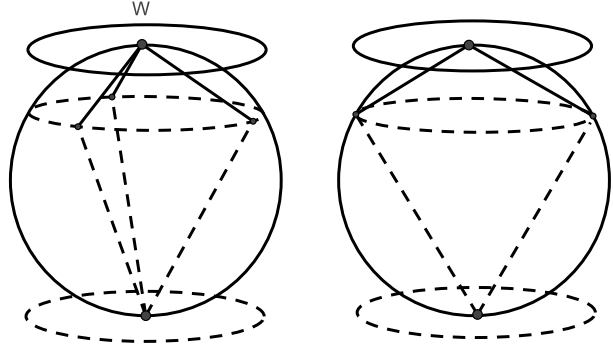}
\caption{Some configurations to satisfy inequality (\ref{minimax4}).}
\end{figure}

\begin{Exm} {\emph{As an example let $G\subset SO(3)$ be the subgroup leaving the north pole $(0,0,1)$ invariant. Since $S_{2}/{G}\ \cong [-1,1]$ the function $f_v^2$ depends only on the last coordinates of the vectors. In this case 
$$f_v^2(w)=2\left(1+\sqrt{(1-v_3^2)(1-w_3^2)}-v_3w_3\right).$$
A simple computation shows that its minimum is taken at $w_3=\pm 1$ (north and south pole) depending on the sign of $v_3$. In particular $m_v^2=2\left(1-v_3\right)$ if $v_3\geq 0$ and $m_v^2=2\left(1+v_3\right)$ if $v_3\leq 0$. Admitting the reflection about the plane $x_3=0$ we can extend the investigation to the group $G\subset O(n)$ leaving the north pole $(0,0,1)$ invariant. In this case 
$$f_v^2(w)=2\left(1+\sqrt{(1-v_3^2)(1-w_3^2)}+\textrm{sgn}(v_3w_3) v_3w_3\right).$$
Its minimum is taken at $w_3=\pm 1$ (north and south pole) if $v_3^2 \leq 1/2$ and $w_3=0$ (equator) if $v_3^2 \geq 1/2$. In particular the minimax values are 
$$m_v^2=2\left(1+\textrm{sgn}(v_3) v_3\right) \qquad \textrm{and} \qquad m_v^2=2\left(1+\sqrt{1-v_3^2}\right),$$ respectively. 
}}
\end{Exm}

\begin{Def}
The group $G\subset O(n)$ has dense orbits if there is an orbit which is a dense subset of the unit sphere. It is called transitive on the Euclidean unit sphere if there is an orbit which covers the entire unit sphere. 
\end{Def}

\begin{Exm}{\emph{The group of rotations about the origin with rational angles has dense orbits in the plane.}}
\end{Exm}

If the subgroup $G$ is transitive on the Euclidean unit sphere then it is impossible to construct a convex body $K$ satisfying (K1) - (K3). Using a continuity argument the same is true for subgroups having dense orbits. It is clear that $G$ has dense orbits if and only if its topological closure $\bar{G}$ is transitive on the unit sphere. In what follows we suppose that $G$ is a {\bf closed} and, consequently, a {\bf compact subgroup} of the orthogonal group. It does not hurt the generality because of $f_v=\bar{f}_v$ (see property (m1)). 

\begin{Lem} Let $w$ be the minimax point of the orbit $P_G(v)$; $f_v^2(w)=4$ if and only if the orbit of $v$ covers the  Euclidean unit sphere, i.e. $G$ is transitive.
\end{Lem}

\begin{Pf} Since $f_v(w)$ is the minimax value, it is clear that for any unit element $z\in S_{n-1}$ 
$$4=f_v^2(w)\leq f_v^2(z)= \sup_{g\in G}\|z-g(v)\|^2\leq 4.$$
This means that $f_v$ is a constant function, i.e. for any $z\in S_{n-1}$ the antipodal point is of the form $g(v)=-z$ and, consequently, the orbit of $v$ covers the entire Euclidean unit sphere.
\end{Pf}

\begin{Lem} Let $w$ be the minimax point of the orbit $P_G(v)$; if $G$ is not transitive on the unit sphere then the set 
\begin{equation}
\label{minimaxset}
M_f(v,w):=\{g(v)\ | \ g\in T_f(v,w)\},\ \ \textrm{where}\ \ T_f(v,w):=\{g\in G\ | \ f_v(w)=\| w-g(v)\| \}
\end{equation}
has at least two different elements or $g(v)=v$ for any $g\in G$.
\end{Lem}

\begin{Pf}
Suppose that we have a uniquely determined element of the form $v_*:=g(v)$, where $g\in T_f(v,w)$. By Corollary 3, $v_*-w$ encloses an angle of measure greater or equal than $90^{\circ}$ with any tangent vector $z$ at $w$. Therefore $w=v_*$ or $w=-v_*$.  In the second case $f_v^2(w)=4$ which is a contradiction because of Lemma 3. Otherwise $w=v_*$ implies that $f_v^2(w)=0$, i.e. $w=g(v)$ for any $g\in G$. In particular $w=v$ under the choice of the identity.  
\end{Pf}

\section{Lazy and busy orbits}

To formulate the geometric meaning of the minimax point process we are going to investigate the converse optimization problem. Let us define the function
$$h_v \colon \mathbb{R}^n\to \mathbb{R},\ \ h_v(w):=\inf_{g\in G} \| w-g(v)\|.$$

\begin{Def}
The maximin point of the orbit $P_G(v):=\{ g(v)\ | \ g\in G\}$ of the unit vector $v$ under the subgroup $G$ is the solution of the optimization problem
\begin{equation}
\label{op3}
\sup h_v\ \ \textrm{subject to}\ \ \|w\|=1;
\end{equation}
if $w$ is a maximin point then 
\begin{equation}
\label{maximinvalue}
r_v:=f_v(w)
\end{equation}
is called the maximin value of $P_G(v)$. 
\end{Def}

It can be easily seen that (\ref{op3}) is equivalent to the optimization problem

\begin{equation}
\label{op4}
\sup h_v^2=\sup \left(\inf_{g\in G} \|w-g(v)\|\right)^2=\sup \inf_{g\in G} \|w-g(v)\|^2  \ \ \textrm{subject to}\ \ \|w\|=1,
\end{equation}
where
$$h_v^2(w)=\inf_{g\in G} \|w-g(v)\|^2.$$
Using Thales theorem
$$\|w-g(v)\|^2+\|w+g(v)\|^2=4$$
and we have that
\begin{equation}
\label{thales}
f_v^2(w)+h_v^2(w)=4\ \ (w\in S_{n-1}), 
\end{equation} 
i.e. if $w$ is a solution of (\ref{op4}) then $-w$ is a solution of (\ref{op2}) and 
$$m_v^2+r_v^2=4.$$
The corresponding monotonicity and stability properties (M1) - (M6) can be also summarized on the model of (m1) - (m6) in Section 2 together with (M7) and (M8):

\begin{itemize}
\item[(M7)] the maximin value (\ref{maximinvalue}) is just the Hausdorff distance of the orbit $P_G(v)$ and the unit sphere $S_{n-1}$,
\item[(M8)] for any unit element $w$
$$h_v^2(w)=\inf_{g\in G} \|w-g(v)\|^2=2\inf_{g\in G} \left(1-\langle w, g(v)\rangle\right)=2\left(1-\sup_{g\in G}\langle w, g(v)\rangle\right),$$
where $\sup_{g\in G}\langle w, g(v)\rangle$ is the so-called support function of the orbit $P_G(v)$ or, in an equivalent way, the support function of its convex hull.
\end{itemize}

Property (M7) motivates the following definition. 

\begin{Def}
The orbit of the element $v\in S_{n-1}$ having maximal Hausdorff distance {\emph{(\ref{maximinvalue})}} from the sphere is called a \emph{lazy orbit}. The orbit of the element $v\in S_{n-1}$ having minimal Hausdorff distance {\emph{(\ref{maximinvalue})}} from the sphere is called a \emph{busy orbit}.
\end{Def}

\begin{Cor}
\label{maximinchar}
Let $w$ be the maximin point of $P_G(v)$; for any tangent vector $z$ to $S_{n-1}$ at $w$
\begin{equation}
\label{minimax5}\max_{g\in T_h(v,w)} \langle g(v)-w, z \rangle \geq 0,\ \ \textrm{where}\ \ T_h(v,w):=\{g\in \bar{G}\ | \ h_v(w)=\| w-g(v)\| \}
\end{equation}
\end{Cor}

Note that $\langle w, z \rangle=0$ because $z$ is tangential to $S_{n-1}$ at $w$. Inequality (\ref{minimax5}) implies that for any $z\in T_{w}S_{n-1}$ there is an element $g\in T_h(v,w)$ such that $z$ and $g(v)-w$ enclose an angle less or equal than $90^{\circ}$. The corresponding property to (m5) says that if $w$ is the maximin point of $P_G(v)$ then $r_w\geq r_v$ and we have the following corollary.

\begin{Cor}
If $w$ is the maximin point of a lazy orbit $P_G(v)$ then $P_G(w)$ is also a lazy orbit with maximin points $g(v)$'s, where $g\in T_h(v,w)$. 
\end{Cor}

\begin{Exc} \label{cone1} {\emph{Using (\ref{thales}) prove that $T_f(v,w)=T_h(v,-w)$. Hint.: if $w$ is the minimax point of $P_G(v)$ then $-w$ is its maximin point and vice versa.}}
\end{Exc} 

 In what follows we suppose that $G$ is a {\bf closed} and, consequently, a {\bf compact subgroup} of the orthogonal group unless otherwise stated. The following results are the analogues of Lemma 3 and Lemma 4, respectively. 

\begin{Lem} Let $w$ be the maximin point of the orbit $P_G(v)$; $h_v^2(w)=0$ if and only if the orbit of $v$ covers the  Euclidean unit sphere, i.e. $G$ is transitive.
\end{Lem}

\begin{Lem} 
\label{atleast2}
Let $w$ be the maximin point of the orbit $P_G(v)$; if $G$ is not transitive on the unit sphere then the set 
\begin{equation}
\label{maximinset}
M_h(v,w):=\{g(v)\ | \ g\in T_h(v,w)\},\ \ \textrm{where}\ \ T_h(v,w):=\{g\in G\ | \ h_v(w)=\| w-g(v)\| \}
\end{equation}
has at least two different elements or $g(v)=v$ for any $g\in G$
\end{Lem}

\begin{Def}
\label{setM}
 Since $M_h(v,w)=M_f(v,w)$ independenty of the choice of the optimization problem we introduce the set
$$M(v,w):=\{g(v)\ | \ g\in T_h(v,w)\}=\{g(v)\ | \ g\in T_f(v,-w)\},$$
where $w$ is the maximin point of the orbit $P_G(v)$ and $G$ is not transitive on the unit sphere.
\end{Def}

In what follows we are going to make the distinction between the reducible and irreducible subgroups in terms of the Hausdorff distance of the orbits from the unit sphere (maximin values).

\begin{Lem}For any bounded $K\subset \mathbb{R}^n$ the convex hull $\textrm{conv\ } \bar{K}$ of the closure of $K$ is the closure of the convex hull of $K$.
\end{Lem} 

\begin{Pf} Since the topological closure $\bar{K}$ of a bounded set $K$ is compact, $K\subset \bar{K}$ implies that  
$$\overline{\textrm{conv\ } K}\subset \textrm{conv\ } \bar{K}$$
because the convex hull of a compact set is compact. If $v\in \textrm{conv\ } \bar{K}$ then we can write 
$$v=\lambda_1v_1+\ldots+\lambda_m v_m$$
as the convex combination of the elements $v_i\in \bar{K}$, i.e. $v_i^k\to v_i$ ($i=1, \ldots, m$) for some sequences of elements in $K$. Therefore 
$$v^k=\lambda_1v_1^k+\ldots+\lambda_m v_m^k$$
is a sequence in $\textrm{conv\ } K$ tending to $v$ and $v\in \overline{\textrm{conv\ } K}$ as was to be proved. 
\end{Pf}

\begin{Rem} {\emph{If $K$ is not bounded then the conclusion is false:  
$$K:= \{ (x,y)\in \mathbb{R}^2\ | \ y=|x|^{-1}\ \}.$$}}
\end{Rem}

\begin{Cor} If $K$ is bounded and $v$ is an extreme point of $\overline{\textrm{conv\ } K}$ then $v\in \bar{K}$.
\end{Cor}

\begin{Pf} If $v\notin \bar{K}$ then $\bar{K}\subset \textrm{conv\ } \bar{K} \setminus \{v\}$ and $\textrm{conv\ } \bar{K} \setminus \{v\}=\overline{\textrm{conv\ } K}\setminus \{v\}$ is a smaller convex set containing $\bar{K}$ than $\textrm{conv\ } \bar{K}=\overline{\textrm{conv\ } K}$. This is a contradiction. 
\end{Pf}

\begin{Lem} 
\label{key0}
A not necessarily closed subgroup $G\subset O(n)$ is irreducible if and only if the origin is an interior point of the topological closure of the convex hulls of non-trivial orbits:
$${\bf 0}\in \ \textrm{int}\ \overline{\textrm{conv\ }\{g(v)\ | \ g\in G\}}$$
for any $v\neq {\bf 0}.$
\end{Lem}

\begin{Pf} We have two cases to discuss:
\begin{itemize}
\item[(i)] ${\bf 0}\notin \overline{\textrm{conv\ }\{g(v)\ | \ g\in G\}},$
\item[(ii)] ${\bf 0}$ is on the boundary of  $\overline{\textrm{conv\ }\{g(v)\ | \ g\in G\}}.$
\end{itemize}
In both cases we are going to construct invariant subspaces of $G$. This means that $G$ is reducible which is a contradiction. Let $v\neq {\bf 0}$ be given and $\bar{K}:=\overline{\textrm{conv\ }\{g(v)\ | \ g\in G\}}.$ 
If ${\bf 0}\notin \bar{K}$ then there exists a uniquely determined closest point to the origin in $\bar{K}$. The unicity implies that it is a fixed point of the transfomations of $G$ because the origin does not move under their actions. In the second case note that ${\bf 0}$ can not be an extreme point of $\bar{K}$ because ${\bf 0}\notin \overline{\{g(v)\ | \ g\in G\}}\subset S_{n-1}$; see Corollary 6. Consider the set $\mathcal{H}$ of the supporting hyperplanes of $\bar{K}$ at the origin. Since the origin is not an extreme point we have a line segment in $\bar{K}$ containing the origin in its relative interior. Such a segment belongs to any $H\in \mathcal{H}$ because ${\bf 0}$ is on the boundary of $\bar{K}$ and $H$ supports $\bar{K}$ at {\bf{0}}. This means that the intersection $L:=\cap_{H\in \mathcal{H}} H$ is at least of dimension $1$. It is invariant under the action of $G$ because for any $g\in G$ we have that
$g(K)=K$ and $g({\bf 0})={\bf 0}$. Therefore the elements of $G$ are varying the elements of $\mathcal{H}$, i.e. their intersection is invariant. The converse of the statement is trivial because the nonempty interior implies the set to be of dimension $n$. Therefore it can not be embedded in a lower dimensional linear subspace.
\end{Pf}

\begin{Thm}
\label{irredhaus}
A not necessarily closed subgroup $G\subset O(n)$ is irreducible if and only if the Hausdorff distance $r_v < \sqrt{2}$ for any $v\in S_{n-1}$.
\end{Thm}

\begin{Pf}
First of all note that $G$ is irreducible if and only if its topological closure is irreducible. Without loss of generality suppose that $G$ is an irreducible, closed and, consequently, compact subgroup. Let $v\in S_{n-1}$ be a given unit element with a maximin point $w$. As we have see above (property (M8))
\begin{equation}
\label{support}
r_v^2=2\left(1-\sup_{g\in G}\langle w, g(v)\rangle\right),
\end{equation}
where $\sup_{g\in G}\langle w, g(v)\rangle$ is the so-called support function of the orbit $P_G(v)$ or, in an equivalent way, the support function of its convex hull. Using Lemma 8, the convex hull contains the origin in its interior, i.e. its support function is the Minkowski functional of the polar body. This means that the support function takes strictly positive values at the nonzero elements of the space. i.e.
$r_v^2 < 2$ as was to be proved. To prove the converse of the statement suppose that $G$ is reducible. i.e. we have an invariant linear subspace $L$. If $v\in L$ and $w\in L^{\bot}\cap S_{n-1}$ then for any $g\in G$
$$\|w-g(v)\|=\sqrt{\|w\|^2+\|g(v)\|^2}=\sqrt{2}\ \ \Rightarrow\ \ h_v(w)= \sqrt{2},$$
i.e. if $G$ is reducible then there exists $v\in S_{n-1}$ such that $r_v\geq \sqrt{2}$. By contraposition and formula (\ref{support}) if $r_v < \sqrt{2}$ for any $v\in S_{n-1}$ then the group is irreducible.
\end{Pf}

\begin{Cor} A not necessarily closed subgroup $G\subset O(n)$ is irreducible if and only if the lazy orbits are closer to the sphere than $\sqrt{2}$. 
\end{Cor}

\begin{Thm} The ratio of the Hausdorff distances $\displaystyle{R:=\sup_{v\in S_{n-1}} r_v}$ and $\displaystyle{r:=\inf_{v\in S_{n-1}} r_v}$ of the lazy and the busy orbits of $G$ to the unit sphere can be estimated as $1/2 \leq r/R \leq 1$. 
\end{Thm}

\begin{Pf} Suppose, in contrary, that $r < R/2$ and let $P_G(v)$ be a lazy orbit with a maximin point $w$, where the maximin value (Hausdorff distance) is attained at: 
$$r_v=\|w- g(v)\|=R\ \ \textrm{for any}\ \ g\in T_h(v,w):=\{g\in G\ | \ h_v(w)=\| w-g(v)\| \}.$$
Choose a mapping $g\in T_h(v,w)$ and consider a busy orbit $P_G(z)$; since the parallel body of $P_G(z)$ with radius $r < R/2$ covers the entire unit sphere there are elements of the form $g_1(z)$ and $g_2(z)$ in the open balls with radius $R/2$ around $g(v)$ and $w$.
Therefore
$$\|w-g_2\circ g_1^{-1}\circ g(v)\|\leq \|w-g_2(z)\|+\|g_2(z)-g_2\circ g_1^{-1}\circ g(v)\|< R/2 + \|z-g_1^{-1}\circ g(v)\|=$$
$$R/2+\|g_1(z)-g(v)\|< R/2+R/2=R,$$
i.e. $g_2\circ g_1^{-1}\circ g(v)$ is a point of the orbit of $v$ which is closer to $w$ than $R$. This is a contradiction; see Figure 2.
\end{Pf}

\begin{figure}
\centering
\includegraphics[scale=0.45]{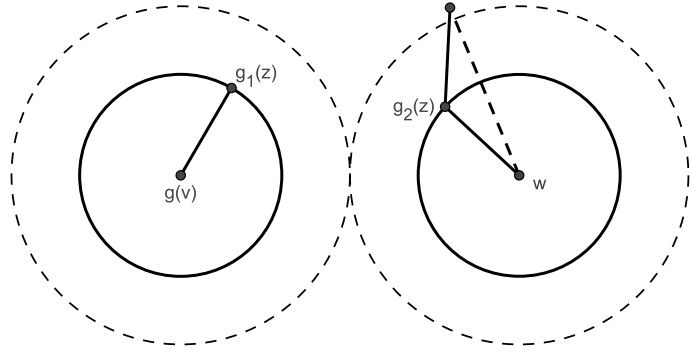}
\caption{The proof of Theorem 3.}
\end{figure}

\subsection{An application of the minimax point process} \cite{VA_2} Let $G\subset O(n)$ be a closed non-transitive irreducible subgroup and consider a unit element $v\in S_{n-1}$. As we have seen above in Lemma 8, its convex hull contains the origin in its interior. To construct a compact convex body containing the origin in its interior such that (K1)-(K3) are satisfied consider the sublevel sets (generalized conics) of the mapping 
\begin{equation}
\label{intfunction}
{\bf x} \mapsto\intop\limits_{\textrm{conv}\ P_G(v)} d_2({\bf x},{\bf y}) \, d{\bf y}.
\end{equation}
Condition (K2) is obviously satisfied because the orbit $P_G(v)$ is invariant under $G$. So is any sublevel set of the mapping (\ref{intfunction}). Unfortunately the explicite computation  of such an integral seems to be impossible in general.  Therefore we follow another way to solve the ellipsoid problem (K1). Let $w$ be the minimax point of $P_G(v)$; recall that $m_v$ denotes the minimax value 
\begin{equation}
\label{minimaxvalue1}
m_v:=\inf_{\|w\|=1}\ f_v(w).
\end{equation}
By the help of the standard calculus \cite{Lee} it can be seen that 
\[s\colon {\bf R}\to {\bf R},\ \ s(t):=\left\{
\begin{array}{rl}
0&\textrm{if} \ \ t\leq m_v\\
(t-m_v)e^{-\frac{1}{t-m_v}}&\textrm{if} \ \ t> m_v
\end{array}
\right.\]
is a smooth convex function on the real line. Let us define the function 
\[u(t):=t+s(t);\]
as we can see nothing happens as far as $t\leq m_v$. If $t>m_v$ then the function $u(t)$ increases its value relative to the argument $t$. Therefore 
$$\intop\limits_{\textrm{conv}\ P_G(v)}  d_2(w,{\bf y}) \, d{\bf y}=\intop\limits_{\textrm{conv}\ P_G(v)} u\circ  d_2(w,{\bf y}) \, d{\bf y}=c_0,$$
i.e. the integrals agree at the minimax point $w$ but at least one of the hypersurfaces
\begin{equation}
\label{conics1}
\intop\limits_{\textrm{conv}\ P_G(v)} d_2({\bf x},{\bf y}) \, d{\bf y}=c_0\ \ \textrm{and}\ \ \intop\limits_{\textrm{conv}\ P_G(v)} u\circ d_2({\bf x},{\bf y}) \, d{\bf y}=c_0
\end{equation}
must be different from the sphere unless the mapping $f_v$ is constant: $f_v(w)=m_v$ for any $w\in S_{n-1}$. Especially $m_v=f_v(-v)=2$ and, by Lemma 3, it is impossible if $G$ is not transitive. Therefore the ellipsoid problem (K1) is solved for irreducible subgroups because invariant ellipsoids under an irreducible subgroup in $O(n)$ must be Euclidean spheres\footnote{Suppose that $G$ contains orthogonal transformations with respect to different inner products. For the sake of simplicity let one of them be the canonical inner product and consider another one given by a symmetric matrix $B$. If $Bx=\lambda x$ for some nonzero vector $x$ then for any vector $y$ and $g\in G$ we have
$$yBg(x)=g^{-1}(y)Bx=\lambda g^{-1}(y)x=\lambda yg(x)\ \ \Rightarrow \ \ Bg(x)=\lambda g(x)$$
because $g^{-1}$ and $g$ are orthogonal transformations with respect to both $B$ and the canonical inner product. Therefore eigenvectors with eigenvalue $\lambda$ form an invariant linear subspace of $G$ which must be the whole space according to the irreducibility. Thus $Bx=\lambda x$ for \emph{all} vectors and the balls with respect to these inner products coincide.}. To satisfy the regularity condition (K3), a sufficiently small increasing of the level rate $c> c_0$ is needed in formula (\ref{conics1}). Then the focal set is contained entirely in the interior of the conic. This means that the regularity condition (K3) is also satisfied. 

\section{An application of the maximin point process: the rank of non-transitive subgroups} 

Let $G\subset O(n)$ be a closed non-transitive subgroup and consider a unit element $v\in S_{n-1}$ with a maximin (or a minimax) point $w$. Recall that the set $M(v,w)$ contains the elements of the orbit $P_G(v)$, where the maximin value (or the minimax value) is attained at; see Definition \ref{setM}. If $w$ is the maximin point then the maximin value is the common distance of the elements in $M(v,w)$ to $w$ and the minimax value is the common distance of the elements in $M(v,w)$ to $-w$.  

\begin{Def} \label{rank}
The subspace spanned by the elements $g(v)\in M(v,w)$ and $w$ is called a flat subspace of $G$ belonging to the orbit of $v\in S_{n-1}$. The rank of the group means the maximal dimension of its flat subspaces, i.e. 
$$\textrm{rank}\ G:=\max_{v\in S_{n-1}} \dim \mathcal{L}\left(M(v,w) \cup \{w\}\right).$$
\end{Def}

\begin{Rem}
{\emph{The terminologies of the flat subspaces and the rank of a group are motivated by Simon's holonomy system theory \cite{Simons} that is the abstract version of the theory of Riemannian holonomies. The main difference comes from its infinitesimal feature due to the associated algebraic curvature tensor to the group that takes the values in the Lie algebra of $G$. The construction is given on the model of the Ambrose-Singer theorem, see also Remark 6. Here we follow a more general and different approach based on convex geometric tools and an optimization problem. The investigations are not restricted by the adjungation of curvature type quadrilinear forms to the group structure.}}
\end{Rem}

\begin{Thm}
\label{redorfinite} If $G\subset O(n)$ is a closed non-transitive subgroup of rank $n$ then it is reducible or finite. 
\end{Thm}

\begin{Pf}
Suppose that $G$ is irreducible and let $v\in S_{n-1}$ be an element, where the rank is attained at. By Theorem \ref{irredhaus} and Corollary \ref{maximinchar} the set $M(v,w)$ does not contain opposite elements and its spherical convex hull 
$$\textrm{conv}_{S_{n-1}} M(v,w):=\textrm{pos} \ M(v,w) \cap S_{n-1}$$
contains the maximin point $w$, where $\textrm{pos}$ means the operator of the positive hull of the sets. The intersection of the images of $\textrm{conv}_{S_{n-1}} M(v,w)$ under different group elements is of measure zero because the spherical convex hull does not contains elements of the form $g(v)$ closer to $w$ than the maximin distance. At the same time the rank is maximal, i.e. $\textrm{conv}_{S_{n-1}} M(v,w)$ is of positive measure on the sphere. Therefore it has only finitely many different images under the action of the elements in $G$. So does the maximin point $w$. By Lemma 8 we can find a basis of the space among the elements of $P_G(w)$ with finitely many possible images. This means that $G$ must be finite. 
\end{Pf}

\begin{Lem} 
\label{lem9}
If $G$ is a proper non-transitive subgroup of the orthogonal group then its rank is at least $2$. 
\end{Lem}

\begin{Pf}
It is a direct consequence of Lemma \ref{atleast2} and the definition of the rank. 
\end{Pf}

\begin{Thm} 
\label{reduciblerank}
If $G\subset O(n)$ is a closed, reducible subgroup then it is finite or its rank is at least $\dim L + 1$, where $L$ is a maximal dimensional invariant subspace such that $G|_{L}$ is irreducible. In particular $\dim L\geq 2$.
\end{Thm}

\begin{Pf} Since $G$ is reducible we can write the space into the direct sum 
$$\mathbb{R}^n=L_1\oplus \ldots \oplus L_d$$
of pairwise orthogonal invariant subspaces. If $\dim L_i=1$ ($i=1, \ldots, d$) then $d=n$ and $G$ must be finite because
$g(v)\equiv (\pm v_1, \ldots, \pm v_n)$ for any $v\in \mathbb{R}^n$, where  $v=v_1+\ldots+v_n$ is the decomposition with respect to the direct sum. Otherwise $\dim L_j \geq 2$ and $G |_{L_j}$ is irreducible for some index $j\in \{1, \ldots, d\}$. Pick a point $v\in L_j$ we have
$$\|w-g(v)\|^2=\|w_1\|^2+\ldots+\|w_j-g(v)\|^2+\ldots+\|w_d\|^2,$$
where $w=w_1+\ldots+w_d$, i.e.
$$\|w-g(v)\|^2=2\left(1-\langle w_j, g(v)\rangle \right)\ \ \Rightarrow\ \ \inf_{g\in G} \|w-g(v)\|^2=2\left(1-\sup_{g\in G}\langle w_j, g(v)\rangle \right)\leq 2$$
because the origin is an interior point of the orbit $P_{G |_{L_j}}(v)$ in $L_j$ due to the irreducibility. Equality occours if and only if $w_j=0$. Therefore $w\equiv (w_1, \ldots, w_{j-1}, 0, w_{j+1}, \ldots, w_d)$ is the solution of the optimization problem
$$\sup \inf_{g\in G} \|w-g(v)\|^2\ \ \textrm{subject to}\ \ \|w\|=1$$
and the flat subspace is spanned by a basis in $L_j$ of the form $g_1(v), \ldots, g_{l_j}(v)$, where $l_j=\dim L_j$ and $w$. 
\end{Pf}

Theorem \ref{redorfinite} and Lemma 9 say that if $n=2$ then any closed non-transitive subgroup is finite. In case of $n=3$ the only possible rank of a not finite closed, irreducible and non-transitive subgroup is $k=2$. The following investigations show that it is impossible.

\begin{figure}
\centering
\includegraphics[scale=0.6]{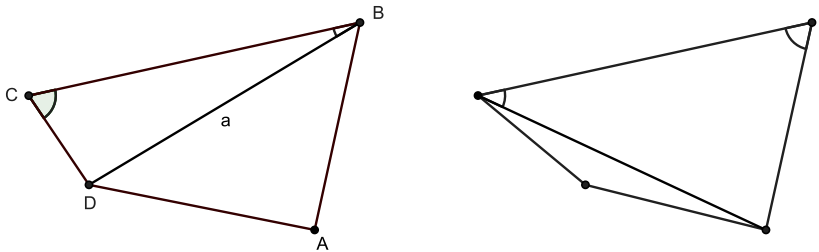}
\caption{Lemma 10.}
\end{figure}

\begin{Lem} Let $ABCD$ be a non-degenerated convex quadrilateral in the plane having equal diagonal segments $AC$ and $BD$. Then each pair of opposite sides contains at least one segment less then the common length of the diagonals.
\end{Lem}

\begin{Pf} Let $a=AC=BD$ be the common length of the diagonals and suppose that $AB\geq a$ and $CD\geq a$. Then, by $CD\geq BD$, we have that $\angle DBC \geq \angle BCD $ (see Figure 3). On the other hand, $AB\geq AC$ implies that $\angle ACB \geq \angle ABC.$ Finally,
$$\angle ABC \leq \angle ACB\leq \angle BCD \leq \angle DBC =\angle ABC-\angle ABD$$
which is a contradiction.
\end{Pf}

The technic of the proof is working in the spherical geometry too.

\begin{Lem}
\label{spherquad} Let $ABCD$ be a non-degenerated convex spherical quadrilateral on the unit sphere $S_2\subset \mathbb{R}^3$ having equal diagonal segments $AC$ and $BD$. Then each pair of opposite sides contains at least one segment less then the common length of the diagonals.
\end{Lem}

\begin{Thm} 
\label{ranktwo}
If $G\subset O(n)$ is a closed, irreducible and non-transitive subgroup of rank two then $g(L)=L$ or $g(L)\cap L=\{\bf{0}\}$ for any $g\in G$, where $L$ is a two-dimensional flat subspace of $G$ belonging to any lazy orbit.
\end{Thm}

\begin{Pf}
By Lemma \ref{atleast2}, the dimension of the flat subspace is exactly two for any $v\in S_{n-1}$ unless $g(v)=v$ ($g\in G$), i.e. the group is reducible. In case of an irreducible subgroup the set $M(v,w)$ contains exactly two different points of the form $z_1=g_1(v)$, $z_2=g_2(v)$ and $w$ (the maximin point) is the midpoint of the geodesic arc between $z_1$ and $z_2$. Suppose that $P_G(v)$ and $P_G(w)$ are lazy orbits such that $z_1$ (or $z_2$) and $w$ are the maximin points of each other (Corollary 5). Therefore the flat subspace $L$ belonging to $w$ is spanned by $z_1$ and $z_2$. In particular, the pairs of the maximin points form a regular polygon inscribed in the intersection $L\cap S_{n-1}$ such that the consecutive vertices belong to the orbits $P_G(v)$ and $P_G(w)$, alternately (see Figure 4).
If $g(L)\neq L$ but $\dim g(L)\cap L=1$ for some element $g\in G$ then there is a spherical quadrilateral with diagonals of equal length and Lemma \ref{spherquad} says that we can choose a side of length less than the common maximin distance of the diagonals such that the consecutive vertices belong to different orbits. This is obviously a contradiction. 
\end{Pf}

\begin{figure}
\centering
\includegraphics[scale=0.45]{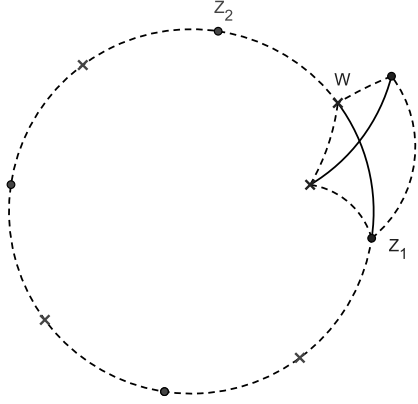}
\caption{The proof of Theorem 6}
\end{figure}

\begin{Rem}{\emph{Using the notations in the proof of the previous theorem we note that the subgroup $G_L\subset G$ leaving the flat subspace $L$ invariant is non-trivial. For example $g:=g_2 \circ g_1^{-1} \in G_L$. Indeed, $g$ sends $z_1$ to $z_2$ and, consequently, $g(w)$ belongs to the linear subspace $L$. Otherwise we need an extra dimension generated by some element of the orbit $P_G(w)$ to keep its maximin point $z_2$ in the spherical convex hull of $M(w, z_2)$. This means that the rank is greater than $2$ which is a contradiction. Therefore we have a one-to-one correspondence $\displaystyle{[g]\in G/G_{L}\mapsto g(L)}$ between the cosets in $G/G_{L}$ and the elements of the orbit $P_G(L)$.}}
\end{Rem}

\begin{Cor}
\label{lowerdim} If $n=2, 3$ then any closed, irreducible and non-transitive subgroup $G\subset O(n)$ is finite.
\end{Cor}

\begin{Rem}
\label{lowerdimrem}
{\emph{The result says that if $n=2, 3$ then any (closed) non-transitive subgroup must be reducible or finite. Using the topological closure it holds for non-transitive subgroups having no dense orbits as well.}}
\end{Rem}

\begin{Cor}
If $n=2, 3$ then there is a finite $G$-invariant system of elements in $\mathbb{R}^n$ for any closed non-transitive subgroup $G\subset O(n)$. 
\end{Cor}

\begin{Pf} If the group is irreducible then Corollary \ref{lowerdim} provides a finite invariant system under $G$ (the orbit of any unit element $v$). Otherwise we always have a one-dimensional invariant subspace by choosing the orthogonal complement of the invariant subspace if necessary, i.e. there exists an invariant system containing the elements $v$, ${\bf 0}$ and $-v$. 
\end{Pf}

\begin{Cor} 
\label{cor10}
If $n=2, 3$ then there is a $G$-invariant polyellipse/polyellipsoid satisfying $(K1)-(K3)$ in $\mathbb{R}^n$ for any closed non-transitive subgroup $G\subset O(n)$.
\end{Cor}

In other words any (closed) non-transitive subgroup $G \subset O(n)$ ($n=2, 3$) can be embedded in the linear isometry group with respect to a non-Euclidean Minkowski functional induced by a polyellipse/polyellipsoid.

\section{An infinitesimal approach: the Lie algebra of non-transitive subgroups}

Let $G\subset O(n)$ be a closed non-transitive subgroup. By the closed subgroup theorem $G$ can be considered as a compact Lie subgroup because of the compactness of the orthogonal group. In what follows we use the following notations: $G^0$ is the unit component, i.e. it is the maximal connected subgroup containing the unit element $e\in G$,
$\displaystyle{d:=\dim G=\dim G^0}$ and  $T_eG$ is the Lie algebra of $G$ considered as the tangent space of $G$ at the identity. As a Lie subalgebra in the space of the skew-symmetric matrices $\mathcal{A}_n(\mathbb{R})$, it is equipped with the usual commutator $[A,B]=AB-BA.$ 

It is well-known that $G^0\subset G$ is compact. On the other hand it is a normal subgroup in $G$ because of the invariance under the action of any inner automorphism -- it follows from the connectedness and the maximality of $G^0$. Finally, the compactness implies that the factor group $G/G^0$ contains only finitely many elements. Now we are going to give a quadratic upper bound for the dimension of $G$ in terms of its rank. Suppose that $v\in S_{n-1}$ is an arbitrarily choosen unit element and let $w$ be the maximin point of the orbit $P_G(v)$. If the elements
$$g_1(v)=v_1, \ldots, g_m(v)=v_m\in M(v,w) \ \ \textrm{and}\ \  v_{m+1}:=w$$
span the flat subspace\footnote{Note that if $G$ is irreducible then the elements $g_1(v)=v_1, \ldots, g_m(v)=v_m$ span the flat subspace because $w$ is in their spherical convex hull; see Corollary 4 and Theorem \ref{irredhaus}.} belonging to $v$ then for any $g\in G$
$$\|g(w)-v_i\|= \|g(w)-g_i(v)\|=\|w-g^{-1}\circ g_i(v)\|\geq \inf_{g\in G}\|w-g(v)\|=h_v(w)=\|w-v_i\|$$
for any $i=1, \ldots, m$, i.e. the elements of the orbit $P_G(w)$ can not decrease the distance from $v_i$'s and the mapping
$$t\mapsto \|e^{tA}(w)-v_i\|^2$$
has zero derivative at $t=0$, where $A\in T_e G$ is a skew-symmetric matrix in the Lie algebra of $G$ and $e^A:=\exp(A)$ for the sake of simplicity. Since
\begin{equation}
\label{expder}
\frac{\partial e^{tA}}{\partial t}=Ae^{tA}=e^{tA}A\ \ \Rightarrow \ \ \frac{\partial e^{tA}}{\partial t}_{t=0} (w)=Aw
\end{equation}
it follows that $\langle Aw, v_i \rangle=0$ for any $i=1, \ldots, m$ and $A\in T_e G$. We have  
\begin{equation}
\label{key1}
\langle Aw, z \rangle \qquad (A\in T_e G, \ z\in L_v),
\end{equation}
i.e. the tangent space of $P_G(w)$ at $w$ is contained in the orthogonal complement of the flat subspace $L_v$ belonging to $v$. 

\begin{Rem}{\emph{Equation (\ref{key1}) motivates to introduce the infinitesimal version of the flat subspace belonging to $v$ as a linear subspace $L$ all of whose elements satisfy (\ref{key1}). Since $L_v\subset L$ it follows that the infinitesimal rank defined as the maximal dimension of the infinitesimal flat subspaces is greater or equal then the rank of the group.}}\end{Rem}

\begin{Thm} For any closed non-transitive subgroup $G\subset O(n)$ 
\begin{equation}
\label{dim1}
\dim G \leq n-k+\frac{(k-1)(k-2)}{2}+\frac{(n-k)(n-k-1)}{2}=k^2-k(n+2)+\frac{n^2+n+2}{2},
\end{equation}
where $k$ is the rank of the group. 
\end{Thm}

\begin{Pf} Let $v$ be a unit element, where the maximal dimension of the flat subspaces is attained at and consider the mapping $\displaystyle{\Phi \colon g\in G \to g(w)\in S_{n-1}}$, where $w$ is the maximin point of $P_G(v)$. Using (\ref{expder}) 
$$\Phi'(e) (A)=Aw$$
and we have, by the rank-nullity theorem for linear mappings, that
\begin{equation}
\label{kernel}
d=\dim \textrm{Ker\ } \Phi'(e)+\dim \textrm{Im\ }\ \Phi'(e) \leq \dim \textrm{Ker} 
\ \Phi'(e)+n-k
\end{equation}
because of equation (\ref{key1}). The Lie subalgebra $\textrm{Ker\ } \ \Phi'(e)$ corresponds to a subgroup $G_k$ in $G$ leaving the element $w$ invariant. This means that the entire flat subspace spanned by $v_1, \ldots, v_m, v_{m+1}$ is also invariant because it is of maximal dimension. Using the restrictions of the elements $G_k$ to the flat subspace $L_v$ of dimension $k$ we have that
$$G_k\subset G_k |_{L_v} \times O(n-k), $$
where, by Theorem \ref{redorfinite}, $G_k |_{L_v}$ is a reducible or a finite subgroup. Therefore its maximal dimension is 
$$\frac{(k-1)(k-2)}{2}=\dim O(k-1),$$
where $k-1$ is the maximal dimension of the non-trivial invariant subspace in $L_v$. Since
$$\dim O(n-k)=\frac{(n-k)(n-k-1)}{2},$$  
the proof is completed.
\end{Pf}

\begin{Rem}
\label{Wang}
{\emph{In the sense of Lemma \ref{lem9}, if $G$ is non-trivial then its rank must be between $2$ and $n$. Taking the right hand side of (\ref{dim1}) as a second order polynomial of $k$ it can be easily seen that}}
\begin{equation}
\label{dim2}
p(2)=p(n)=\frac{(n-1)(n-2)}{2}, \ \ \textrm{{\emph{where}}}\ \ p(k):=k^2-k(n+2)+\frac{n^2+n+2}{2}.
\end{equation} 
{\emph{Therefore the common value of $p(2)$ and $p(n)$ is the maximal dimension of a closed non-transitive subgroup in $O(n)$. The bound also appears in Wang's theorem \cite{Wang} up to the additive term of the free parameters for the translation part of the isometries.}}
\end{Rem}

\begin{Lem}
\label{pointwise}
If there is a pointwise fixed linear subspace under the action of the unit component $G^0$ then $G$ is a finite or a reducible group.
\end{Lem}

\begin{Pf}
Suppose that $u$ is a unit element such that $g(u)=u$ for any $g\in G^0$. Since $G/G^0$ is finite $u$ also has only finitely many images under the action of $G$ because $g(u)$ is independent of the representation of the equivalence class for any $g\in G$. If $G$ is irreducible then we have a basis among the finitely element of $P_G(u)$ (see Lemma \ref{key0}) and any $g\in G$ sends each element of the bases into the finite set $P_G(u)$. Therefore $G$ is finite. By contraposition, if $G$ is not finite then it must be reducible.  
\end{Pf}

\begin{Thm}
\label{ranknminusone}
If a closed non-transitive subgroup $G\subset O(n)$ is of rank $n-1$ then we have the following possible cases:
\begin{itemize}
\item[(i)] it is a finite or a reducible group,
\item[(ii)] $2d=n$, where $d\geq 2$ and there exists a decomposition $L_1\oplus \cdots \oplus L_d$ of the space into the direct sum of two-dimensional linear subspaces such that the unit component $G^0$ acts transitively on each two-dimensional subspace $L_i$, where $i=1, \ldots, d$. In particular $G^0$ is a compact connected Abelian subgroup, 
$2\leq \left |G/G^0 \right |\leq d!2^d$ and $g(L_i)=L_j$ for some indices $i$ and $j$ independently of the representation of the elements in $G/G^0$.
\end{itemize}
\end{Thm}

\begin{Pf} Suppose that $G$ is irreducible and let $w$ be the maximin point of the orbit $P_G(v)$ of a unit element $v\in S_{n-1}$, where the maximal dimension of the flat subspaces is attained at. If $u$ is the unit normal of the flat subspace $L_v$ of dimension $n-1$ it follows, by (\ref{key1}) that $Aw=\langle Aw,u\rangle u$ for any $A\in T_eG$. Therefore
$$[A,B]w=\langle Bw,u\rangle Au-\langle Aw, u\rangle Bu, \ \ \textrm{i.e.}\ \ [A,B]w \ \bot \ u\ \ \textrm{and}\ \ [A,B]w \ \| \ u$$
at the same time. This means that $[A,B]w=0$. Repeating the process at $g(w)$ ($g\in G$) as the maximin point of $P_G(v)$ we have, by Lemma \ref{key0}, that $[A,B]=0$. Therefore $G^0$ is a compact connected Abelian subgroup. If $G_0$ is trivial then $G$ must be finite by a compactness argument. Otherwise we have a non-zero $A\in T_eG$. Taking the canonical form of the orthogonal transormations suppose that $L$ is a two-dimensional invariant subspace of $g_0=e^A$ such that $g|_{L}$ acts transitively on the subspace $L$ because of $A\neq 0$. Since $[A,B]=0$ is equivalent to the commutativity $e^Ae^B=e^Be^A$ it follows that the elements of $G^0$ sends $L$ into a two-dimensional invariant subspace of $g_0$. Joining $g_0$ and $g\in G^0$ with a continuous curve\footnote{Since $G^0$ is a differentiable manifold, the connectedness and the arcwise connectedness are equivalent.} $c\colon [0,1]\to G^0$ we have that $c(t)(L)$ is a continuous mapping into the Grassmannian manifold taking at most finitely many values, i.e. it is constant:
$$L=c(0)(L)=c(1)(L)=g(L).$$ 
This means that $L$ is an invariant two-dimensional subspace for any $g\in G^0$. Since the dimension is $d\geq 1$ we can repeat the same argument $d$ times to present a decomposition of the space into the direct sum of pairwise orthogonal invariant linear subspaces
$$\mathbb{R}^n=L_0\oplus L_1\oplus \ldots \oplus L_d,$$
where $G^0$ acts on each two-dimensional subspace $L_i$ ($i=1, \ldots, d$) transitively and $L_0$ is a pointwise fixed linear subspace under the action of $G^0$. If it is non-trivial then, by Lemma \ref{pointwise}, we have that $G$ is a finite or a reducible group. Otherwise $n=2d$ and we can reduce the decomposition of the space to the form 
\begin{equation}
\label{reduceddecomposition}
\mathbb{R}^n=L_1\oplus \ldots \oplus L_d.
\end{equation}
If $d=1$ then $n=2$ and $G$ must be finite or reducible in the sense of Corollary 8. Otherwise $d\geq 2$. Since each inner automorphism by an element of $G$ sends $G^0$ into a connected subgroup containing the identity, it follows that $G^0$ is invariant under the action of the inner automorphisms. Therefore $g(L_i)=L_j$ for some indices $i$ and $j$ independently of the representations of the elements in $G/G^0$. We must have at least two different equivalence classes to avoid the reducible case $G=G^0$. On the other hand there is a homomorphism 
$$P\colon [g]\in G/G^0 \mapsto P[g]:=(g(L_1), \ldots, g(L_d))$$
of $G/G^0$ into the symmetry group of order $d$. The factorization by the kernel of $P$ identifies two different classes in $G/G^0$ if and only if they are represented by the elements $g_1$ and $g_2$ such that
$$g_1(L_i)=g_2(L_i) \qquad (i=1, \ldots, d)$$
and $g_2\circ g_1^{-1}$ is an orientation-reversing mapping of at least one of the invariant subspaces $L_1, \ldots, L_d$. Recall that $G^0$ acts transitively on each two-dimensional invariant subspace. Therefore the possible values of any equivalence class identified with $[e]$ by the kernel of $P$ is $(\pm L_1, \ldots, \pm L_d),$ where the sign refer to the possible orientations of the invariant linear subspaces. 
\end{Pf}

\begin{Rem}
\label{abelianflat}{\emph{Since the mapping $A\in T_e G \mapsto \langle Aw,u\rangle $
is a linear functional, it follows that we have a ($d-1$) - dimensional Lie subalgebra $T_eH \subset T_e G$ such that
$Aw=0$ for any $A\in T_e H$, i.e. there is a (compact) connected subgroup $H^0 \subset G^0$ leaving the maximin point $w$ invariant. Taking 
\begin{equation}
\label{component}
w=w_1+\ldots+w_d, \quad w_i\in L_i \quad (i=1, \ldots, d)
\end{equation}
it is clear that $w_i=0$ for any indices $i=1, \ldots, d$ except exactly one component to provide a subgroup of dimension $d-1$ leaving $w$ invariant. Suppose that (for example) $w\in L_d$. If the group $G$ is irreducible then there are linearly independent elements 
$$g_1(v), \ldots, g_{n-1}(v)$$
to span the flat subspace $L_v$ of dimension $n-1$ and $w$ is in their spherical convex hull. Then the orthogonal projection of each $g_i(v)$ into $L_d$ must be of the form $\lambda w$ for some $\lambda\geq 0$. Indeed, suppose that  the orthogonal projection of $g_j(v)$ into $L_d$ is not in the positive hull of $w$ for some index $j=1, \ldots, n-1$. Using a non-trivial rotation $g'|_{L_d}$ of the projected vector in $L_d$ such that the ray of the rotated element and $w$ coincide, the distance of $w$ and $P_G(v)$ is succesfully decreased:
$$\|w-g'\circ g_j(v)\|< \|w-g_j(v)\|$$
due to the Pythagorean theorem. This is a contradiction  because the right hand side is the minimal distance of $w$ and $P_G(v)$. We have that the flat subspace $L_v$ is the direct sum of $L_1, \ldots, L_{d-1}$ and the one-dimensional linear subspace generated by $w\in L_d$.}}
\end{Rem}

\begin{figure}
\centering
\includegraphics[scale=0.65]{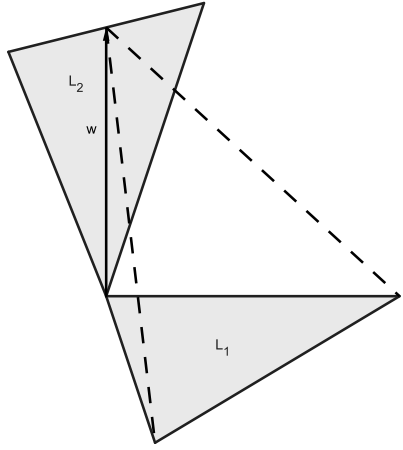}
\caption{The space is of dimension $4$; see Remark \ref{abelianflat}}
\end{figure}

\begin{Cor}
Suppose that $n \geq 3$ and let $G\subset O(n)$ be a closed non-transitive subgroup such that $1\leq d < n/2$ or $d> n/2$. If $G$  is of rank $n-1$ then it is reducible. 
\end{Cor}

\begin{Cor}
If $n$ is an odd number and $G\subset O(n)$ is a closed non-transitive subgroup of rank $n-1$ then it is finite or reducible. 
\end{Cor}

\subsection{Non-transitive subgroups in the four-dimensional case} According to Remark \ref{Wang} the maximal dimension of a closed non-transitive subgroup $G\subset O(4)$ is $3$. If $d=0$ then $G^0$ is trivial and $G$ must be finite. In case of $d=1$ we have that $G^0$ must be a (compact) Abelian group and we can repeat the argument of the proof of Theorem \ref{ranknminusone} to conclude that there exists a decomposition $L_1\oplus L_2$ of the space into the direct sum of two-dimensional linear subspaces such that the elements of $G^0$ are of the form
\begin{equation}
\label{2dimcon}
\left( {\begin{array}{cc}
   R_{\omega_1} (t)&0\\ 
   0& R_{\omega_2}(t)\\
  \end{array} } \right), \quad \textrm{where} \quad R_{\omega_i}(t)=\left( {\begin{array}{cc}
\ \ \ \cos (\omega_i t)& -\sin (\omega_i t)\\ 
   -\sin (\omega_i t)& \ \ \ \cos (\omega_i t)\\
  \end{array} } \right) \quad (i=1,2).
	\end{equation}
On the other hand 
\begin{multline}
\label{2dimdiscon}
\begin{aligned}
\left( {\begin{array}{cccc}
   1&0&0&0\\ 
	   0&1&0&0\\ 
		   0&0&1&0\\ 
			   0&0&0&1\\ 
  \end{array} } \right), \left( {\begin{array}{cccc}
   1&0&0&0\\ 
	   0&-1&0&0\\ 
		   0&0&1&0\\ 
			   0&0&0&1\\ 
  \end{array} } \right), \left( {\begin{array}{cccc}
   1&0&0&0\\ 
	   0&1&0&0\\ 
		   0&0&1&0\\ 
			   0&0&0&-1\\ 
  \end{array} } \right), \left( {\begin{array}{cccc}
   1&0&0&0\\ 
	   0&-1&0&0\\ 
		   0&0&1&0\\ 
			   0&0&0&-1\\ 
  \end{array} } \right)\\
	\\
	\left( {\begin{array}{cccc}
   0&0&1&0\\ 
	   0&0&0&1\\ 
		   1&0&0&0\\ 
			   0&1&0&0\\ 
  \end{array} } \right), \left( {\begin{array}{cccc}
   0&0&1&0\\ 
	   0&0&0&-1\\ 
		   1&0&0&0\\ 
			   0&1&0&0\\ 
  \end{array} } \right), \left( {\begin{array}{cccc}
   0&0&1&0\\ 
	   0&0&0&1\\ 
		   1&0&0&0\\ 
			   0&-1&0&0\\ 
  \end{array} } \right), \left( {\begin{array}{cccc}
   0&0&1&0\\ 
	   0&0&0&-1\\ 
		   1&0&0&0\\ 
			   0&-1&0&0\\ 
  \end{array} } \right)
	\end{aligned}
	\end{multline}
	represent the possible classes in $G/G^0$. If $G$ is not reducible then $k=2$ (Corollary 11, Theorem 4) and $\omega_1\neq 0$, $\omega_2\neq 0$ (Lemma \ref{pointwise}). To discuss the case of $d=2$ we need the following inductive argument.
	
\subsection{An inductive characterization} Suppose that $G$ is a closed irreducible subgroup of dimension $1\leq d \leq n-1$. Choosing a unit element $v\in S_{n-1}$ let us introduce the inner product
\begin{equation}
\label{averageinner}
i(A,B):=\int_{\textrm{conv}\ P_G(v)} \langle Au, Bu\rangle\, du
\end{equation}
for the Lie algebra $T_eG$. It is obviously positive definite because of Lemma \ref{key0}. On the other hand it is invariant under the action $I_g(A):=gAg^{-1}$ of the inner automorphism by any element $g\in G$, i.e.
$i(I_g(A),I_g(B))=i(A,B).$ This means that the mapping $I\colon g\in G^0 \mapsto I_g$ is a homomorphism of $G^0$ into the orthogonal group of a Euclidean space of dimension $d\leq n-1$. On the other hand $I(G^0)$ is a compact connected subgroup in $SO(d)$. The first isomorphism thereom implies that
$\textrm{Ker}\ I$ is a normal subgroup in $G^0$. On the other hand the factor group $G^0/\textrm{Ker}\ I$ is isometric to a compact connected subgroup in $SO(d)$ as the special orthogonal group of $T_eG$ equipped with the inner product (\ref{averageinner}). If it is transitive then $G^0/\textrm{Ker}\ I$ can be found among the elements of the table in Section 1 that contains the unit components of transitive subgroups. Otherwise  we can use the characterization of non-transitive subgroups in dimension $1, \ldots, n-1$. 

\subsection{The case of $n=4$ and $d=2$} We restrict our investigations to the case of irreducible subgroups (see  Theorem \ref{reduciblerank} and Theorem \ref{ranknminusone}). First of all note that $\dim G^0=2$ implies that $\textrm{Ker}\ I$ is at least of dimension one because of $\dim SO(2)=1$. This means that we have a non-zero element $B$ in $T_eG$ such that $gAg^{-1}=A \quad (A\in T_eG),$
where $g=e^B$ (note that the exponential map is surjective onto the identity component in case of compact Lie groups). In other words the one-parameter subgroup $ge^{tA}g^{-1}$ is generated by $A$, i.e. $g$ commutes with the elements of $G^0$. Therefore $[A,B]=0$ for any $A\in T_eG$. This means that $G^0$ is a compact Abelian group because it is of dimension two. The same argument as in the proof of Theorem \ref{ranknminusone} gives the decomposition $\mathbb{R}^4=L_1\oplus L_2$ of the space into the direct sum of two-dimensional invariant linear subspaces under the action of $G^0$. Therefore it is generated by the one-parameter subgroups
\begin{equation}
\label{2dimconagain}
\left( {\begin{array}{cc}
   R_{\omega_1} (t)&0\\ 
   0& R_{\omega_2}(t)\\
  \end{array} } \right), \ \left( {\begin{array}{cc}
   R_{\omega_3} (s)&0\\ 
   0& R_{\omega_4}(s)\\
  \end{array} } \right)
	\end{equation}
	where
	$$ R_{\omega_i}(t)=\left( {\begin{array}{cc}
\ \ \ \cos (\omega_i t)& -\sin (\omega_i t)\\ 
   -\sin (\omega_i t)& \ \ \ \cos (\omega_i t)\\
  \end{array} } \right) \quad (i=1,2),$$
		$$R_{\omega_i}(s)=\left( {\begin{array}{cc}
\ \ \ \cos (\omega_i s)& -\sin (\omega_i s)\\ 
   -\sin (\omega_i s)& \ \ \ \cos (\omega_i s)\\
  \end{array} } \right) \quad (i=3,4)
	$$
	and $(\omega_1, \omega_2)$, $(\omega_3, \omega_4)$ are linearly independent pairs; (\ref{2dimdiscon}) represents the possible classes in $G/G^0$. If $G=G^0$ then we have a reducible group of rank at least $3$ (Theorem 5). Therefore we can suppose that there is a mapping $g\in G$ exchanging the subspaces $L_1$ and $L_2$. Consider a unit element $v=v_1+v_2$ ($v_1\in L_1, v_2\in L_2$). Since we can use independent rotations in $L_1$ and $L_2$ we have
	$$h_v^2(w)= \min \left \{ \big| \|w_1\|-\|v_1\|\big |^2+\big | \|w_2\|-\|v_2\|\big |^2, \big | \|w_1\|-\|v_2\|\big |^2+\big | \|w_2\|-\|v_1\|\big |^2 \right \}.$$
	Using the angle representations
	$$\|v_1\|=\cos \alpha_0,\ \|v_2\|=\sin \alpha_0, \|w_1\|=\cos \alpha, \|w_2\|=\sin \alpha,$$
 it follows that we should solve the optimization problems
	\begin{equation}
	\label{opt1}
	\begin{aligned}
	\max_{\alpha} \ 2-2\cos \alpha_0 \cos \alpha -2\sin \alpha_0 \sin \alpha\ \ \textrm{subject to}\\
	\ 2-2\cos \alpha_0 \cos \alpha -2\sin \alpha_0 \sin \alpha\leq 2-2\cos \alpha_0 \sin \alpha -2\sin \alpha_0 \cos \alpha,
	\end{aligned}
	\end{equation}
		where $\alpha_0, \alpha\in [0, \pi/2]$ and
			\begin{equation}
	\label{opt2}
	\begin{aligned}
	\max_{\alpha}\  2-2\cos \alpha_0 \sin \alpha -2\sin \alpha_0 \cos \alpha\ \ \textrm{subject to}\\
	2-2\cos \alpha_0 \sin \alpha -2\sin \alpha_0 \cos \alpha\leq 2-2\cos \alpha_0 \cos \alpha -2\sin \alpha_0 \sin \alpha.
	\end{aligned}
	\end{equation}
 Reformulating the problems
	\begin{equation}
	\label{opt1new}
	\max_{\alpha} \ 1-\cos (\alpha_0-\alpha) \ \textrm{subject to}\ 	\ 0\leq (\cos \alpha_0 - \sin \alpha_0)(\cos \alpha-\sin \alpha),
	\end{equation}
		where $\alpha_0, \alpha\in [0, \pi/2]$ and
			\begin{equation}
	\label{opt2new}
	\max_{\alpha}\  1-\sin (\alpha_0+\alpha)\ \ \textrm{subject to} \ \ 	 0\geq (\cos \alpha_0 - \sin \alpha_0)(\cos \alpha-\sin \alpha).
	\end{equation}
	The following table shows the possible cases of the solutions.
\vspace{0.3cm}
\noindent
\begin{center}
\begin{tabular}{|l|l|l|l|}
\hline
 &&&\\
$\alpha_0$ & problem (\ref{opt1new})& problem (\ref{opt2new})& the square of the maximin dist.\\
\hline
 &&&\\
$0\leq\alpha_0 \leq \pi/8$ &$\alpha=\pi/4$&$\alpha=\pi/4$&$2\left(1-\frac{\ 1}{\sqrt{2}}\cos \alpha_0 -\frac{\ 1}{\sqrt{2}}\sin \alpha_0\right)$\\
\hline
 &&&\\
$\pi/8 \leq \alpha_0\leq \pi/4$&$\alpha=0$&$\alpha=\pi/2$&$2(1-\cos \alpha_0)$\\
\hline
 &&&\\
$\pi/4 \leq \alpha_0 \leq 3\pi/8$&$\alpha=\pi/2$&$\alpha=0$&$2(1-\sin \alpha_0)$\\
\hline
 &&&\\
$3\pi/8 \leq \alpha_0 \leq \pi/2$&$\alpha=\pi/4$&$\alpha=\pi/4$&$2\left(1-\frac{\ 1}{\sqrt{2}}\cos \alpha_0 -\frac{\ 1}{\sqrt{2}}\sin \alpha_0\right)$\\
\hline
\end{tabular}
\end{center}
\vspace{0.3cm}
The busy orbits belong to $\alpha_0=\pi/8$ or $3\pi/8$ with Hausdorff distance
$$\sqrt{2\left(1-\frac{\cos \alpha_0 +\sin \alpha_0}{\sqrt{2}}\right)}\approx 0.3901806441.$$	
The lazy orbits occurs at $\alpha_0=0, \pi/4$ or $\pi/2$ with Hausdorff distance 
$$\sqrt{2-\sqrt{2}}\approx 0.7653668650.$$
Suppose that $k=2$, $\alpha_0=0$, i.e. (for example) $v=(1,0,0,0)\in L_1$ and the flat subspace belonging to the lazy orbit is spanned by
$z_1=g_1(v)$, $z_2=g_2(v)$ and $w$, where the maximin point $w$ bisects the geodesic segment from $z_1\in L_1$ to $z_2\in L_2$ (for example). Since we can use independent rotations in $L_1$ and $L_2$, it is possible to send $z_1$ to its original position up to $\pm$ such that $z_2$ does not arrive to $\pm z_2$. It contradicts to the synchronicity property Theorem \ref{ranktwo}. Therefore the rank must be $3$. 

\subsection{The case of $n=4$ and $d=3$} Let $v$ be a unit vector such that $P_G(v)$ is a lazy orbit and suppose that $w$ is its maximin point. If $k=2$ and $b_1:=w$, $b_2, \ldots, b_4$ is an orthonormal basis such that the flat subspace $L_v$ belonging to $v$ is spanned by $b_1$ and $b_2$ then, for any $A\in T_eG$ we have by (\ref{key1}) that 
$$\langle Ab_1, b_i\rangle =0 \quad (i=1,2).$$
Therefore the matrix representation of any element of the Lie algebra is of the form
\begin{equation}
\label{toranktwo}
\left( {\begin{array}{cccc}
   0&0&a&b\\ 
	   0&0&c&d\\ 
		   -a&-c&0&s\\ 
			   -b&-d&-s&0\\ 
  \end{array} } \right)
	\end{equation}
Since $\dim G=3$, (\ref{kernel}) implies that we have an at least one-dimensional subgroup in $G^0$ leaving $b_1=w$ invariant. In case of $k=2$ such a subgroup can not move the flat subspace because of Theorem \ref{ranktwo}, i.e. $b_2$ is also a fixed element. Therefore 
\begin{equation}
\left( {\begin{array}{cccc}
   0&0&0&0\\ 
	   0&0&0&0\\ 
		   0&0&0&1\\ 
			   0&0&-1&0\\ 
  \end{array} } \right)\in T_e G.
	\end{equation}
	Some direct computations give the special form of matrices spanning the Lie algebra of $G$:
	$$E_1:= \left( \begin {array}{cccc} 0&0&0&0\\ \noalign{\medskip}0&0&0&0
\\ \noalign{\medskip}0&0&0&1\\ \noalign{\medskip}0&0&-1&0\end {array}
 \right), \ E_2:= \left( \begin {array}{cccc} 0&0&x&y\\ \noalign{\medskip}0&0&z&w
\\ \noalign{\medskip}-x&-z&0&s\\ \noalign{\medskip}-y&-w&-s&0
\end {array} \right), $$
\begin{equation}
\label{exc}
E_3:=[E_1, E_2]= \left( \begin {array}{cccc} 0&0&y&-x\\ \noalign{\medskip}0&0&w&-z
\\ \noalign{\medskip}-y&-w&0&0\\ \noalign{\medskip}x&z&0&0\end {array}
 \right)
\end{equation}
and, consequently,
$$[E_1,E_3]= \left( \begin {array}{cccc} 0&0&-x&-y\\ \noalign{\medskip}0&0&-z&-w
\\ \noalign{\medskip}x&z&0&0\\ \noalign{\medskip}y&w&0&0\end {array}
 \right)=-E_2+sE_1,
$$
$$[E_2, E_3]= \left( \begin {array}{cccc} 0&-2\,xw+2\,yz&-sx&-sy
\\ \noalign{\medskip}2\,xw-2\,yz&0&-sz&-sw\\ \noalign{\medskip}sx&sz&0
&{w}^{2}+{x}^{2}+{y}^{2}+{z}^{2}\\ \noalign{\medskip}sy&sw&-{w}^{2}-{x
}^{2}-{y}^{2}-{z}^{2}&0\end {array} \right).
$$
This means that $E_1, E_2$ and $E_3$ generate a Lie algebra if and only if $z=tx$, $w=ty$ for some common proportional term $t\in \mathbb{R}$ provided that $x$ and $y$ do not vanish simultaneously:
\begin{equation}
\label{ranktwodimthree}
\begin{aligned}
 M_1&=\left( \begin {array}{cccc} 0&0&0&0\\ \noalign{\medskip}0&0&0&0
\\ \noalign{\medskip}0&0&0&1\\ \noalign{\medskip}0&0&-1&0\end {array}
 \right),  \ M_2=\left( \begin {array}{cccc} 0&0&x&y\\ \noalign{\medskip}0&0&tx&ty
\\ \noalign{\medskip}-x&-tx&0&s\\ \noalign{\medskip}-y&-ty&-s&0
\end {array} \right),\\
 M_3&:=[M_1, M_2]=\left( \begin {array}{cccc} 0&0&y&-x\\ \noalign{\medskip}0&0&ty&-tx
\\ \noalign{\medskip}-y&-ty&0&0\\ \noalign{\medskip}x&tx&0&0
\end {array} \right),\\
[M_1, M_3]&=-M_2+sM_1, \ [M_2, M_3]=-sM_2+(x^2+y^2)(t^2+1)M_1.
\end{aligned}
	\end{equation}
Another possible case is $x=y=0$ in (\ref{exc}), i.e. $w$ is a fixed vector under the action of $G^0$. In the sense of Lemma \ref{pointwise} it follows that $G$ is reducible and $k\geq 3$ (Theorem 5). 

\begin{Def}
The group $G$ is called locally reducible/irreducible if its unit component $G^0$ is a reducible/irreducible group. 
\end{Def}

\begin{Rem}
{\emph{It is clear that if $G$ is reducible, then it is locally reducible. By contraposition a locally irreducible group is irreducible.}}
\end{Rem}

\subsection{Summary} The following table summarizes the possible cases of closed nontransitive groups in case of $n=4$. 
\vspace{0.3cm}
\noindent
\begin{center}
\begin{tabular}{|l|l|l|l|l|}
\hline
 \multicolumn{5}{|c|}{}\\
 \multicolumn{5}{|c|}{Non-transitive closed sobgroups in case of $n=4$} \\
\multicolumn{5}{|c|}{}
\\
\hline
&&&&\\
dimension &$d=0$&$d=1$&$d=2$&$d=3$\\
rank &&&&\\
\hline
&&&&\\
$k=2$&finite& locally reducible:&--& locally irreducible:\\
&&(\ref{2dimcon}) and (\ref{2dimdiscon})&& (\ref{ranktwodimthree})\\
\hline
&&&&\\
$k=3$&finite&reducible: Thm. \ref{ranknminusone}&locally reducible: & reducible: Thm. \ref{ranknminusone}\\
&&&(\ref{2dimconagain}) and (\ref{2dimdiscon})&\\
\hline
&&&&\\
$k=4$&finite&reducible: Thm. \ref{redorfinite}&reducible: Thm. \ref{redorfinite}&reducible: Thm. \ref{redorfinite}\\
\hline
\end{tabular}
\end{center}

\section{Concluding remarks} Let $M$ be a connected Riemannian manifold and consider a (not necessarily torsion-free) metric linear connection $\nabla$ on $M$; $G:=H_p$ denotes the holonomy group of $\nabla$ at the point $p$. The unit component is denoted by $H_p^0$. The rank of the connection $\nabla$ is defined as the rank of the topological closure of $H_p$. The connectedness provides that the rank is independent of the choice of $p\in M$ because the parallel transports preserve the flat subspaces. Therefore each result of the previous sections can be interpreted in the context of metric linear connections. If $H_p$ has no dense orbits on the Euclidean unit sphere in $T_pM$ then we can construct an $H_p$ - invariant generalized conic body \cite{VA_2}. The $H_p$ - invariance allows us to extend it by parallel transports to any point of the manifold; see Z. I. Szab\'{o}'s idea \cite{Szab1}. The smoothly varying family of compact convex bodies provides a non-Riemannian metric environment for $\nabla$: the Minkowski functionals induced by the generalized conics in the tangent spaces constitute a so-called Finslerian fundamental function such that the parallel transports with respect to $\nabla$ preserve the Finslerian length of tangent vectors (compatibility property). The Finsler geometry is the alternative of the Riemanninan geometry for the linear connection $\nabla$. The result can be formulated as follows.  

\begin{Thm} \emph{\cite{VA_2}}
Suppose that $M$ is a connected Riemannian manifold and $\nabla$ is a metric linear connection on M. If $p\in M$ and the holonomy group $H_p$ of $\nabla$ has no dense orbits on the Euclidean unit sphere in $T_p M$ then there is a non-Riemannian Finsler manifold equipped with the fundamental function $F\colon TM\to \mathbb{R}$ such that the parallel transports with respect to $\nabla$ preserve the Finslerian length of tangent vectors and the indicatrix bodies in the tangent spaces are generalized conics.
\end{Thm}

\begin{Def}
Finsler manifolds admitting compatible linear connections are called generalized Berwald manifolds. If the compatible linear connection has zero torsion then we have a classical Berwald manifold.
\end{Def}

\begin{Rem}
\label{nonbalanced}
{\emph{For a Riemannian manifold the indicatrix hypersurfaces are conics (quadratic hypersurfaces) in the classical sense. In the sense of the previous theorem, if $M$ is a non-Riemannian generalized Berwald manifold then the indicatrix bodies can be supposed to be generalized conics instead of the classical ones.}} 
\end{Rem}

Using Corollary \ref{lowerdim} we have the following statement in case of lower dimensional generalized Berwald manifolds; the corresponding statements for (classical) Berwald manifolds can be found in \cite{Szab1}. 

\begin{Cor} The compatible linear connection has zero curvature for any two-dimensional Finsler manifold.  
In case of a three-dimensional Finsler manifold any compatible linear connection is reducible or it has zero curvature. 
\end{Cor}

Some further reformulations of the results in the previous sections:
\begin{itemize} 
\item If the compatible linear connection is of maximal rank $n=\dim M$ then it is reducible or it has zero curvature (see Theorem 4 and the Ambrose-Singer theorem about the holonomy of linear connections).
\item If the compatible linear connection is of rank $n-1$ then it is locally reducible or it has zero curvature (Theorem 8). 
\item If $\dim M$ is an odd number and the compatible linear connection is of rank $n-1$ then it is reducible or it has zero curvature (Corollary 12). 
\item In case of a four-dimensional Finsler manifold any compatible linear connection is locally reducible  or it has zero curvature except the case of the Lie algebra (\ref{ranktwodimthree}) of its holonomy group. 
\end{itemize}

\begin{Rem}
The so-called simply connectedness is a standard topological condition for the manifold to eliminate the difference between the locally reducible and the reducible linear connections. 
\end{Rem}

\subsection{An open problem} Construct a linear connection on a four-dimensional connected Riemannian manifold with Lie algebra (\ref{ranktwodimthree}) of its holonomy group.

\end{document}